\documentclass[11pt,a4paper,reqno]{amsart}
\usepackage{graphicx} 

\usepackage{amsmath}
\usepackage[foot]{amsaddr}
\usepackage{enumitem}
\usepackage{xcolor}
\usepackage{amssymb}
\usepackage{amsthm}
\usepackage{array}
\newcolumntype{M}[1]{>{\centering\arraybackslash}m{#1}}

\usepackage{float}

\usepackage{cuted}
\usepackage{color}
\usepackage{graphicx}

\usepackage{tikz}
\usetikzlibrary{shapes.geometric}
\usetikzlibrary{arrows.meta,arrows}

\usepackage{tabularx}

\usepackage{comment}

\usepackage[hidelinks, colorlinks=false]{hyperref}

\usepackage{cleveref}

\usepackage[normalem]{ulem}

\newtheorem{thm}{Theorem}
\newtheorem{remark}[thm]{Remark}

\newtheorem{definition}[thm]{Definition}

\newcommand{\R}{\mathbb{R}}
\newcommand{\N}{\mathbb{N}}

\newcommand{\n}[1]{\left\Vert #1\right\Vert}

\newcommand{\eqb}{\begin{flalign}}
\newcommand{\eqe}{\end{flalign}}
\newcommand{\Ra}{\Rightarrow}
\newcommand{\Lra}{\Leftrightarrow}

\allowdisplaybreaks

\definecolor{mygrey}{RGB}{220, 220, 220}
\definecolor{concrete}{RGB}{136, 150, 136}
\definecolor{steel}{RGB}{200, 200, 200}
\definecolor{copper}{RGB}{220, 127, 100}
\definecolor{redd}{RGB}{200, 50, 50}
\definecolor{yellow}{RGB}{255, 254, 235}
\definecolor{green}{RGB}{213, 243, 231}
\definecolor{green2}{RGB}{169, 209, 142}
\definecolor{lightgrey}{RGB}{240, 240, 240}
\definecolor{lightred}{RGB}{255, 244, 239}
\definecolor{lightlilac}{RGB}{239, 239, 255}
\definecolor{wood}{RGB}{133, 94, 66}
\definecolor{BackgroundGrey}{RGB}{210, 210, 210}

\begin{document}
\title[Spatial decay of perturbations in hyperbolic equations]{Spatial decay of perturbations in hyperbolic equations with optimal boundary control}
\thanks{B.O.\ extends his gratitude to the Thüringer Graduiertenförderung for funding his scholarship.}

\author{Benedikt Oppeneiger$^1$} 
\author{Manuel Schaller$^2$} 
\author{Karl Worthmann$^1$}

\address{$^1$Optimization-based Control Group, Institute of Mathematics, Technische Universität Ilmenau (e-mail: $\{$benedikt-florian.oppeneiger, karl.worthmann$\}$@tu-ilmenau.de).}
\address{$^2$Faculty of Mathematics, Chemnitz University of Technology (e-mail: manuel.schaller@math.tu-chemnitz.de).}

\begin{abstract}                
Recently, 
domain-uniform stabilizability and detectability has been the central assumption 
to ensure robustness in the sense of exponential decay of spatially localized perturbations in optimally controlled evolution equations. In the present paper we analyze a chain of transport equations with boundary and point controls with regard to this property.
Both for Dirichlet and Neumann boundary and coupling conditions, we show a necessary and sufficient criterion on control domains which allow for the domain-uniform stabilization of this equation. We illustrate the results by means of a numerical example.
\end{abstract}


\maketitle
\section{Introduction}

Robustness with regard to perturbations is a key issue in finding practicable solutions for the design of controllers which can be applied to real-world systems. In recent years robustness of optimal trajectories has therefore been a major research topic in the context of optimal control for partial differential equations. A manifestation of robustness with respect to perturbation of initial values is the turnpike property, which states that optimal trajectories stay near an optimal steady state for most of the time, see e.g.~\cite{Damm2014, TrelZuaz15} for nonlinear problems,
\cite{Breiten2020,Schaller2020} for PDEs or the overview articles ~\cite{Faulwasser2022Turnpike,GeshZuaz22}.
In~\cite{Schaller2020} it was shown that stabilizability and detectability imply localized sensitivities in time w.r.t.~ perturbations of the dynamics: temporally localized source terms only have a local effect on the optimal solution.

In \cite{Shin2022}, spatially localized robustness was analyzed for the case of graph-structured NLPs using first- and second-order optimality conditions uniform in the size of the graph. PDEs were considered only recently in~\cite{Goettlich2025,Goettlich2024}, where the assumption of domain-uniform stabilizability (and its dual notion of detectability) has been para\-mount in proving that spatially localized perturbations only cause spatially localized deviations in the behaviour of optimally controlled hyperbolic equations. 
Loosely speaking,  and considering an evolution equation on $L^2(\Omega)$, $\Omega \subset \R^n$ with $n\in \N$, domain-uniform stabilizability means that for any initial value there is a control such that the resulting state satisfies 
\begin{align}\label{eq:stab}
    \|x(t)\|_{L^2(\Omega)}\leq Me^{-\omega t}\|x^0\|_{L^2(\Omega)}
\end{align}
for constants $M \geq 1 $ and $\omega > 0$ uniform in the size of the spatial domain $\Omega$. Here, the uniformity of the constants in the domain size is crucial to obtain the exponential spatial decay in optimal solutions by a scaling argument.

Domain-uniform stabilizability and detectability was characterized for distributed control of transport equations in ~\cite{Goettlich2025}.
In this paper we extend these results to the case of boundary and point control. We show that the transport equation with boundary control is domain uniformly stabilizable if and only if the distance between two neighboring control access points is bounded from above uniformly in the domain size.

\noindent The paper is structured as follows: In Section~\ref{Sec: ProbStat} we introduce the system under consideration and the notion of domain-uniform stabilizability. In Section~\ref{Sec: Dirichlet} we derive a state space representation of the system with Dirichlet boundary control and proof a characterization of domain-uniform stabilizability for this case. In Section~\ref{Sec: Neumann} we extend these results to the case of Neumann boundary control. Finally we illustrate our results and their implications on optimally controlled transport equations via some numerical simulations.

\section{Problem Statement}
\label{Sec: ProbStat}
\begin{figure}[htb]
    \centering
    \includegraphics[width=\linewidth]{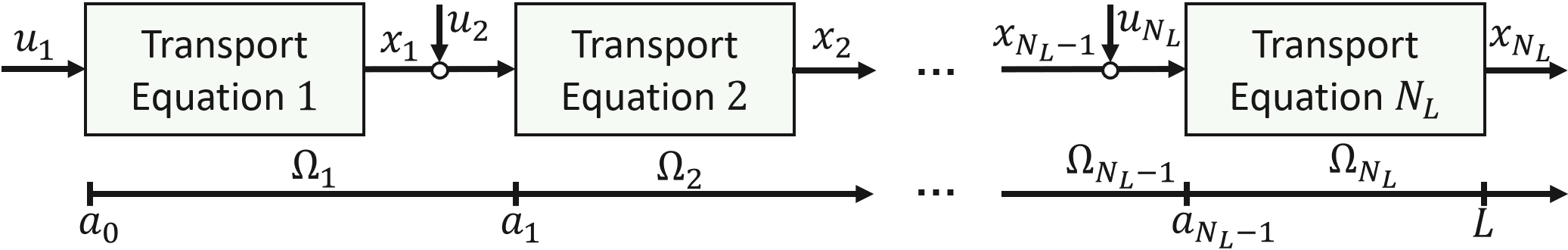}
    \caption{Sketch of a chain of coupled transport equations}
    \label{fig: Sketch}
\end{figure}

In this paper we analyze a chain of transport equations as illustrated in Figure~\ref{fig: Sketch} w.r.t.
 domain-uniform stabilizability and detectability. 

Similar to, e.g.~\cite{Xu2023}, the terminology \emph{chain} refers to a directed tree with only one path. As depicted in Figure~\ref{fig: Sketch}, consecutive transport equations are coupled in such a manner that the successor's left boundary value depends on the corresponding control input and its predecessor's right boundary value. The boundary condition of the leftmost transport equation only depends on the input since it has no predecessor. 

In this work, we provide a criterion on the positioning of the control access points between the individual subsystems such that the coupled system is domain-uniformly stabilizable in the sense of \eqref{eq:stab}.

\noindent To make the setting depicted in Figure~\ref{fig: Sketch} mathematically precise, we define the following: 
Let $(a_i)_{i\in \N_0} \subset \mathbb{R}$ be a strictly increasing sequence with 
$a_0 = 0$ and $\Delta := \inf_{i\in \N} \,l_i > 0$, $l_i = | \Omega_i |$ with $\Omega_i :=(a_{i-1},a_i)$. 
For $L>0$, we define $N_L := \min_{i\in \N} \{i\in \N: a_i\geq L\}$. 
On $\Omega ^L := (0,L)$
and for $i\in \{1,\ldots ,N_L\}$, we consider the transport equations
\begin{subequations}
\label{eq: Transport}
\begin{align}
\label{eq: TransportState}
    \forall (\omega,t) \in \left(\Omega_i\cap \Omega ^L\right) \times \R _{\geq 0}: \dot{x}_i(\omega,t)&=-cx_i'(\omega,t)\\
    \label{eq: Initial}
    \forall \, \omega \in \Omega_i\cap \Omega ^L:\, x_i(\omega,0) &= x_L^0 |_{\Omega _i \cap \Omega ^L},
\end{align}
\end{subequations}

with initial distribution $x_L^0\in L^2(\Omega ^L)$, constant velocity $c > 0$ and $x_0(0,t):=0 \ \forall \, t \in \R _{\geq 0}$.

In the subsequent analysis we characterize domain-uniform stabilizability for the two cases of Dirichlet coupling conditions
\begin{equation}
    \label{eq: DirichletBC}
    x_i(a _{i-1},t) = x_{i-1}(a _{i-1},t) + u_i(t)
\end{equation}
and Neumann coupling conditions
\begin{equation}
    \label{eq: NeumannBC}
    \frac{\partial}{\partial \omega} x_i(a _{i-1},t) = \frac{\partial}{\partial \omega} x_{i-1}(a _{i-1},t) + u_i(t).
\end{equation}

where we have a control $u \in L^2(\R _{\geq 0},\R^{N_L})$.
Let $X_{L,i}\subset L^2(\Omega _i\cap \Omega ^L)$ be a Hilbert space such that the differential operator 
\begin{equation}
\label{eq: AOperator}
    A_{L,i}: D(A_{L,i})\subset X_{L,i} \rightarrow X_{L,i}, \quad x \mapsto -c\frac{\partial}{\partial \omega}x.
\end{equation}
generates a strongly continuous semigroup $\left( T_{L,i}(t)\right)_{t\geq 0}$ on $X_{L,i}$ (e.g.\ $X_{L,i}:=L^2(\Omega ^L\cap \Omega _i)$ for Dirichlet and $X_{L,i}:=H^1(\Omega ^L\cap \Omega _i)$ for Neumann boundary conditions).
    We call $x_{i}^u = x_{i}^u(\cdot, \cdot; x_L^0) \in C(\R _{\geq 0}, X_{L,i})$ a mild solution of the $i$-th transport equation in~\eqref{eq: Transport} if and only if
\begin{equation*}
    x_i^u(\cdot,t)-x_i^u(\cdot,0) = A_{L,i} \int _0^t x_i^u(\cdot,\tau)\,\mathrm{d}\tau \quad \forall\,t \geq 0.
\end{equation*}

By $X_L$ we define the space of functions $v \in L^2(\Omega ^L)$ which fulfill $v|_{\Omega _i \cap \Omega ^L} \in X_{L,i}$ for all $i\in \{1,\ldots , N_L\}$. By using the inner product
\begin{equation*}
    \langle v,w\rangle _{X_L} = \sum _{i=1}^{N_L} \langle v|_{\Omega _i \cap \Omega ^L},w|_{\Omega _i \cap \Omega ^L}\rangle _{X_{L,i}}
\end{equation*}
it is easy to see, that $X_L$ is a Hilbert space.
The overall solution of the chain is always denoted by $x_u\in C(\R _{\geq 0}, X_L)$, where
$(x_u(t))(\omega):= (x_i^u(t))(\omega)$, $\omega \in (a_{i-1},a_i]$.
For readability, we mostly write $x(\omega,t)$ instead of $x(t)(\omega)$. 

In the following we define the notion of domain-uniform stabilizability adapted from~\cite[Def.~9]{Goettlich2025}.
\begin{definition}[Domain-uniform stabilizability]
   \label{Def: StabDect}
   For $\gamma \geq 0$ we call the transport equation in~\eqref{eq: Transport} $\gamma$-domain uniformly exponentially stabilizable if and only if there exist constants $M,k>0$ such that for all $x^0_L\in X_L$ there is $u_L \in L^2(\R _{\geq 0},\R^{N_L})$
   such that the solution $x_{u_L}(t,\omega)=x_{u_L}(t,\omega; x_L^0)$ fulfills the inequality
    \begin{equation*}
        \forall  L>\gamma \, \forall x_0^L \in Z_L\, \forall  t \geq 0: \n{x_{u_L}(\cdot,t)}_{X_L} \leq Me^{-kt}\n{x_L^0}_{X_L}.
    \end{equation*}
\end{definition}
Our goal is to characterize those sequences $(a_i)_{i \in \N _0}$ for which the system defined in~\eqref{eq: Transport} and~\eqref{eq: DirichletBC} respectively~\eqref{eq: NeumannBC} is 0-domain-uniformly exponentially stabilizable.

\section{Dirichlet boundary control}
\label{Sec: Dirichlet}

In this section, we 
derive a necessary and sufficient criterion for the domain-uniform stabilizability of the transport equation with Dirichlet boundary control. For this purpose we 
first rewrite our system in a state-space form 
to derive a mild solution formula. Then, 
we present our main result, i.e., a characterization of domain-uniform stabilizability.

\noindent\textbf{State space form and solution formula}. 
We consider the initial value problem~\eqref{eq: Transport} with Dirichlet boundary condition~\eqref{eq: DirichletBC} with state space $X_{L,i}:=L^2(\Omega _i \cap \Omega ^L)$ and $v_i(t):=x_{i-1}(a_{i-1},t)+u_i(t)$. 
We choose $D(A_{L,i}):= \{x \in H^1(\Omega _i \cap \Omega ^L): x(a_{i-1})=0\}$ as domain of the differential operator~$A_{L,i}$ in~\eqref{eq: AOperator}.
In the following, we often leave out the subscripts $L,i$ for readability. 
The operator~$A$ is the infinitesimal generator of the strongly continuous semigroup $\left(T(t)\right)_{t\geq 0}$ on $X$ given by
\begin{equation}
\label{eq: DirichletSemigroup}
    \left(T(t)y\right)(\omega)=\left\{ \begin{array}{cc}
        0, & \omega \leq ct \\
        y(\omega-ct),    & \mathrm{else}
    \end{array}
    \right..  
\end{equation}

Therefore the mild solution of the autonomous ($v_i\equiv 0$) equations is given by~\cite[p. 61]{Jacob2012} \begin{equation*}
    x_{i}^0: \R _{\geq 0} \rightarrow X_{L,i}, \,\, x_{i,x_L^0}^0(\omega,t)
    =\left\{ \begin{array}{cc}
        0, & \omega \leq ct \\
        x_L^0(\omega-ct),    & \mathrm{else}
        \end{array}
    \right..  
\end{equation*}

For $\beta \in \sigma (A)$ we denote by $X_{-1}$  the completion of the space $X$ with respect to the norm $\| ( \beta I - A )^{-1} \cdot \|_X$ (which indeed is independent of $\beta$)
and by $X_1^d$ the normed space $\left({D}(A^*),\n{(\beta^*I-A^*)\cdot}_X\right)$. 
Using $X_1^d\subset H^1(\Omega _i \cap \Omega ^L) \hookrightarrow C(\overline{\Omega _i \cap \Omega ^L})$ we define the delta distribution $\delta _\alpha \in X_{-1}$ via 
\begin{equation*}
    \forall \psi \in X_1^d: \left\langle \delta _\alpha,\psi \right\rangle_{X_{-1},X_1^d} = \psi (\alpha) \qquad\forall\,\alpha \in \Omega _i \cap \Omega ^L.
\end{equation*}
Following~\cite[Section 10.1]{Weiss2009}, we rewrite the transport equation with Dirichlet boundary condition in the state-space form
\begin{equation}
\label{eq: CauchyDirichlet}
    \dot{x_i} = A_{L,i,-1}x_i+B_{L,i}v_i, \quad x_i(0)=x_L^0|_{\Omega _i \cap \Omega ^L},
\end{equation}
where $A_{L,i,-1}$ is the unique extension of $A_{L,i}$ to $X_{L,i,-1}$ and the input operators are given by
\begin{equation}
\label{eq: InputOperator}
    B_{L,i}: \R \rightarrow X_{-1,L,i} \quad\text{ and }\quad  B_{L,i}u:=u\delta_{a_{i-1}}.
\end{equation}
Using this state space representation in combination with the variation of constants formula it is possible to derive the mild solution of \eqref{eq: CauchyDirichlet}:
\begin{equation}
    \label{eq: DirichletSolution}
    x_i^{u_i}(\omega ,t) = \begin{cases}
        w_i(\omega,t), & \omega \!\leq\! a_{i-1}\!+\!ct \\
        x_L^0(a_{i-1}+\omega - ct),    & \mathrm{else}
    \end{cases}
\end{equation}
where
\begin{equation*}
    w_i(\omega,t):=x_{i-1}^{u_{i-1}}\!\left(a_{i-1},t\!-\!\frac{\omega -a_{i-1}}{c}\right)\!+\!u_i\!\left(t\!-\!\frac{\omega-a_{i-1}}{c}\right).
\end{equation*}
The overall solution $x^u$ can be computed via the formula
\begin{equation*}
    x_i^{u_i}(\omega ,t)\!=\! \chi _{[ct,L]} x_0^L(\omega - ct) \sum _{i=1}^{N_L} \chi _{I_i(t)}(\omega)u_i\left(t\!-\!\frac{\omega-a_{i-1}}{c}\right)
\end{equation*}
where $\chi _\Omega$ is the characteristic function of a set $\Omega \subset \R$ and
\begin{equation}
\label{eq: IntervalDef}
    I_i(t):=\Omega ^L \cap (a_{i-1},a_{i-1}+ct].
\end{equation}

\noindent\textbf{Characterization of domain-uniform stabilizability}.
The mild solution formula in~\eqref{eq: DirichletSolution} shows that for a given point $t$ in time the control $u$ only influences the solution of the chain of transport equations in~\eqref{eq: Transport} on the set
\begin{equation*}
    U_t :=  \Omega ^L \cap \left(\overset{N_L}{\underset{{i=1}}{\cup}} I_i(t)\right) = \Omega ^L \cap \left(\overset{N_L}{\underset{{i=1}}{\cup}}(a_{i-1},a_{i-1}+ct]\right).
\end{equation*}
Therefore, for any control $u \in L^2(\R _{\geq 0},\R^{N_L})$, the time required to cause an exponential decay of the initial value function everywhere on the spatial domain, increases as the maxi\-mum distance between two neighboring elements of the sequence $\left(a_i\right)_{i\in \N}$ becomes larger.
Using this observation we can proof a necessary and sufficient criterion on the control domain $\Omega _c:= \cup_{i\in \N _0}a_i$ such system~\eqref{eq: Transport} with boundary condition~\eqref{eq: DirichletBC} becomes domain-uniformly stabilizable.
\begin{thm}
\label{Thm: DirichletStab}
    Assume $x_0^L \in H^1(0,L)$. Then the following statements are equivalent:
    \begin{itemize}
        \item [(i)] The chain of transport equations~\eqref{eq: Transport} with Dirichlet boundary control~\eqref{eq: DirichletBC} is $0$-domain-uniformly exponentially stabilizable.
        \item [(ii)]  The sequence $(a_i)_{i\in \N}$ satisfies the condition
        \begin{equation*}
            \exists \, L_0 > 0 \, \forall i \in \N _0: a_{i+1}-a_i \leq L_0.
        \end{equation*}
        \item [(iii)] There exists $L_0>0$ such that for all intervals $I\subset \R _{\geq 0}$
        \begin{equation*}
            I \cap \Omega _c = \emptyset \implies |I| \leq L_0.
        \end{equation*}
    \end{itemize}
\end{thm}
\begin{proof}
    We first show the equivalence (i) $\Leftrightarrow$ (ii).\\
    \textbf{(i) $\Rightarrow$ (ii):} Assume that (ii) is not fulfilled, i.e. we find $\epsilon >0$ such that
    \begin{equation*}
        \forall \, L_0>0 \, \exists \, i_{L_0} \in \N _0 : a_{i_{L_0}+1} - a_{i_{L_0}} > L_0+\varepsilon.
    \end{equation*}
    For arbitrary $M,k>0$ choose
    \begin{equation*}
        L_0 = 3c \left|\frac{\ln(M)}{k}\right| \, \mathrm{and} \, t_{L_0}= \frac{L_0}{2c}.
    \end{equation*}
    Using this definition we find
    \begin{equation}
    \label{eq: DirichletExpEstimate}
        M e^{-kt_{L_0}}=Me^{-1.5k\left|\frac{\ln(M)}{k}\right|}<Me^{-k\frac{\ln(M)}{k}}=1.
    \end{equation}
    Now choose $L>a_{i_{L_0}+1}$ and $x_0 \in X_L\setminus \{0\}$ such that
    \begin{equation}
        \label{eq: Dirichletx0Choice}
        \chi _{\left[a_{i_{L_0}}, a_{i_{L_0}}+\varepsilon\right]} x_0 = x_0.
    \end{equation}
    With this specific choice we find for an arbitrary control $u_L \in L^2(\R _{\geq 0},\R^{N_L})$ the inequalities
    \begin{flalign*}
        &\n{x_u(\cdot,t_{L_0})}_{X_L} \geq \n{\chi _{[a_{i_{L_0}}+ct_{L_0}, a_{i_{L_0}}+ct_{L_0}+\varepsilon]}x_{x_0}^u(\cdot,t_{L_0})}_{X_L}\\
        = &\n{x_0}_{L^2(a_{i_{L_0}}, a_{i_{L_0}}+\varepsilon)}
        \overset{\eqref{eq: Dirichletx0Choice}}{=} \n{x_0}_{X_L}\overset{\eqref{eq: DirichletExpEstimate}}{>} M e^{-kt_{L_0}}\n{x_0}_{X_L}
    \end{flalign*}
    where the second equality follows from the solution formula in~\eqref{eq: DirichletSolution}. Therefore the chain of transport equations~\eqref{eq: Transport} with Dirichlet boundary control~\eqref{eq: DirichletBC} is \textit{not} domain-uniformly stabilizable.\\
    \textbf{(ii) $\Rightarrow$ (i):} 
    First note that for a control $u_L\in C(\R _{\geq 0},\R^{N_L})$ which fulfills the conditions
\begin{equation*}
    u_{L,1}(0)=x_0(0)  \quad \textrm{and} \quad \forall i\in \{2,\ldots , N_L\}: u_{L,i}(0)=0
\end{equation*}
equation~\eqref{eq: DirichletSolution} implies
\begin{equation*}
    \forall i\in \{1,\ldots , N_L\}:x_i^{u_L}(t,\cdot) \in C(\Omega ^L).
\end{equation*}

Define $\Delta t := \frac{\Delta}{c}$ and
\begin{equation}
\label{eq: Feedback}
    \kappa: \R_{\geq 0} \rightarrow \R, \quad \kappa(t) := \begin{cases}
        2t/\Delta t, & t \leq \Delta t/2 \\
        1, & \mathrm{else}
    \end{cases}
\end{equation}
    For $L>0$ we define a state feedback $u_L\in C(\R _{\geq 0},\R^{N_L})$,
    \begin{equation*}
        u_1(t):=(1-\kappa (t))x_0(ct), \quad u_i(t):=\kappa(t)x_{i-1}(a_{i-1},t).
    \end{equation*}
    Using the formula in~\eqref{eq: DirichletSolution} we find
    \begin{equation*}
        x_i^{u_L}(\cdot ,t)\equiv 0 \qquad\forall\, t \geq 2 L_0/c \geq L_0/c + \Delta t / 2.
    \end{equation*}
    Furthermore, $\| x_{u_L}(t) \|_{X_L} 
        \leq \| x_L^0 \|_{X_L}$ holds for $t \leq 2 L_0/c$.
    Therefore, we find for any $k>0$ and $M:=e^{2k\frac{L_0}{c}}$ that
    \begin{equation*}
        \forall \, L>0 \, \forall \, t \geq 0: \n{x_{u_L}(\cdot,t)}_{X_{\Omega ^L}} \leq Me^{-kt}\n{x_L^0}_{X_{\Omega ^L}}.
    \end{equation*}
    This shows the $0$-domain-uniform stabilizability of~\eqref{eq: Transport}.\\

    \textbf{(ii) $\Rightarrow$ (iii):} Assume that (ii) is fulfilled. Let $I \subset \R_{\geq 0}$ an interval such that $I\cap \Omega = \emptyset$. Since $(a_i)_{i\in \N}$ is an unbounded sequence, $I$ is bounded and we can define
    \begin{equation*}
        \underline{i}_I:=\max \{i: a_i \leq \inf I\} \quad \textrm{and} \quad \overline{i}_I:=\min \{i: a_i \geq \sup I\}.
    \end{equation*}
    Since $I\cap \Omega = \emptyset$ we have $\overline{i}_I=\underline{i}_I+1$. Therefore we find
    \begin{equation*}
        |I|=\sup I - \inf I\leq a_{\underline{i}_I+1}-a_{\underline{i}_I} \leq L_0.
    \end{equation*}
    \textbf{(iii) $\Rightarrow$ (ii):} Assume that condition (ii) is not fulfilled, i.e. for every $L_0 >0$ we find $i_{L_0}\in \N$ such that
    \begin{equation*}
        a_{i_{L_0}}-a_{i_{L_0}-1}> L_0.
    \end{equation*}
    Define $I_{L_0}:=(a_{i_{L_0}-1},a_{i_{L_0}})$. Then we have $I\cap \Omega_c = \emptyset$ and $|I|>L_0$. Therefore condition (iii) is also not satisfied.\hfill $\Box$
\end{proof}

\begin{remark}
\label{Detectability}
A system is domain-uniformly detectable if and only if the dual system is stabilizable.
For $\mathrm{D}(A_{L,i}^*):=\{x \in H^1(\Omega _i \cap \Omega _L): x(a_i)=0\}$ the dual of~\eqref{eq: AOperator} is given by
\begin{equation*}
        A_{L,i}^*: \mathrm{D}(A_{L,i}^*)\subset X_{L,i} \rightarrow X_{L,i}, \quad A_{L,i}^*x\mapsto c\frac{\partial}{\partial \omega}x.
\end{equation*}
The dual of~\eqref{eq: InputOperator} is the observation operator
\begin{equation*}
        C_{L,i}: X_{1,L,i}^d \rightarrow \R, \quad C_{L,i}\Psi:=\Psi (a_{i-1}).
\end{equation*}
Therefore using this observation operator domain-uniform detectability of the transport equation is equivalent to domain-uniform stabilizability of a transport equation with reversed direction of transport. Since the direction of transport is irrelevant for  domain-uniform stabilizability Theorem~\eqref{Thm: DirichletStab} can also be used to characterize domain-uniform detectability.
\end{remark}

\section{Neumann boundary control}
\label{Sec: Neumann}
In the main result of this section we will characterize the control domains $\Omega _c=\underset{i\in \N}{\cup}a_i$ which allow for the domain-uniform stabilization of the transport equation with Neumann boundary control. For this purpose we first transform the system with Neumann boundary control into a system with weakly differentiable Dirichlet boundary control. Using a similar procedure as in Section~\ref{Sec: Dirichlet} we are then able to derive a suitable solution formula which allows us to show our result on domain-uniform-stabilizability.

\noindent\textbf{State space form and solution formula}.
In this section we consider the initial value problem~\eqref{eq: Transport} with Neumann boundary condition~\eqref{eq: NeumannBC}. We choose the state space $X_{L,i}:=H^1(\Omega _i \cap \Omega ^L)$ and $v_i(t):=\frac{\partial}{\partial \omega}x_{i-1}(a_{i-1},t)+u_i(t)$ where we only allow for control inputs $u_i\in H^1(\R_{\geq 0},\R)$.

By integrating the boundary condition~\eqref{eq: NeumannBC} in time we find the corresponding Dirichlet condition since
\begin{flalign}
\label{eq: Neumann2Dirichlet}
\begin{split}
    x_i(a_{i-1},t) = \int _0^t \frac{\partial}{\partial \tau} x_i(a_{i-1},\tau)\mathrm{d}\tau + x_i(a_{i-1},0)
    &\overset{\eqref{eq: TransportState}}{=} \int _0^t -c\frac{\partial}{\partial \omega} x_i(a_{i-1},\tau)\mathrm{d}\tau + x_i(a_{i-1},0)\\
    &\underset{\eqref{eq: NeumannBC}}{\overset{\eqref{eq: Initial}}{=}} -c\int _0^t v_i(\tau)\mathrm{d}\tau + x_L^0(a_{i-1}).
\end{split}
\end{flalign}
Using the input transformation
\begin{equation}
\label{eq: InputTrafo}
    \mathfrak{v}_i(t):=-c\int _0^t v_i(\tau)\mathrm{d}\tau + x_L^0(a_{i-1})
\end{equation}
we can therefore rewrite~\eqref{eq: NeumannBC} as a Dirichlet boundary condition
\begin{equation}
\label{eq: DirichletTransformed}
    x_i(a_{i-1},t) = \mathfrak{v}_i(t).
\end{equation}
Note that the space of feasible transformed inputs $\mathfrak{v}_i$ is given by $H^2(\R_{\geq 0},\R)$. The domain of the differential ope\-rator $A_{L,i}$ from~\eqref{eq: AOperator} corresponding to~\eqref{eq: DirichletTransformed} is given by the set $D(A_{L,i}):= \left\{x \in H^2(\Omega _i \cap \Omega ^L): \frac{\partial}{\partial \omega}x(a_{i-1})=0\right\}$.

Using a similar argument as in Section~\ref{Sec: Dirichlet} we find the mild solution
\begin{equation*}
    x_i^0: \R _{\geq 0} \rightarrow X _{L,i}, \, x_i^0(\omega,t)
    =\left\{ \begin{array}{cc}
        x_L^0(a_{i-1}), & \omega \leq ct \\
        x_L^0(\omega-ct),    & \mathrm{else}
        \end{array}
    \right..  
\end{equation*}
of the autonomous equation ($v_i\equiv 0$).
Again following the approach in~\cite[Section 10.1]{Weiss2009} the transport equation with Dirichlet boundary condition~\eqref{eq: DirichletTransformed} can be rewritten into the state space form~\eqref{eq: CauchyDirichlet} with input operator~\eqref{eq: InputOperator}.
This leads to the mild solution formula
\begin{equation}
    \label{eq: NeumannSolution}
    x_i^{u_i}(\omega ,t) =\left\{\! \begin{array}{cc}
        x_L^0(a_{i - 1}) - c\int _0^{t - \frac{\omega  -  a_{i - 1}}{c}} v_i(\tau)\mathrm{d}\tau, & \omega \!\leq\! a_{i -  1}+ct \\
        x_L^0(a_{i - 1} + \omega - ct),    & \mathrm{else}
    \end{array}
    \right. 
\end{equation}
for $i\in \{1,\ldots , N_L\}$. The overall solution $x_u$ is given by
\begin{flalign*}
    x_i^{u_i}(\omega ,t) &= \chi _{[ct,L]} x_0^L(\omega - ct)
    \quad + \sum _{i=1}^{N_L} \chi _{I_i(t)}(\omega)\int _0^{t -  \frac{\omega -  a_{i -  1}}{c}} u_i\left(\tau \right) \mathrm{d}\tau
\end{flalign*}
where $I_i(t)$ is defined as in~\eqref{eq: IntervalDef}. 

\noindent\textbf{Characterization of domain-uniform stabilizability}.
Using the mild solution formula~\eqref{eq: NeumannSolution} we can prove a similar result as Theorem~\ref{Thm: DirichletStab} for the Neumann case.
\begin{thm}
    Assume $x_0^L \in H^2([0,L])$. Then the following statements are equivalent:
    \begin{itemize}
        \item [(i)] The chain of transport equations~\eqref{eq: Transport} with Neumann boundary control~\eqref{eq: NeumannBC} is domain-uniformly stabilizable.
        \item [(ii)] The sequence $(a_i)_{i\in \N}$ satisfies the condition
        \begin{equation*}
            \exists \, L_0 > 0 \, \forall i \in \N _0: a_{i+1}-a_i \leq L_0.
        \end{equation*}
        \item [(iii)] There exists $L_0>0$ such that for all intervals $I\subset \R _{\geq 0}$
        \begin{equation*}
            I \cap \Omega _c = \emptyset \implies |I| \leq L_0.
        \end{equation*}
    \end{itemize}
\end{thm}
\begin{proof} 
We show (ii) $\Ra$ (i) since (i) $\Ra$ (ii) $\Lra$ (iii) directly follows from the arguments used in the proof of Theorem~\ref{Thm: DirichletStab}.\\
\textbf{(ii) $\Ra$ (i):}
Define $\Delta t := \frac{\Delta}{c}$ and $\kappa: \R_{\geq 0} \rightarrow \R$ as in~\eqref{eq: Feedback}.
Note that for $r\geq 0$
\begin{flalign*}
    \int _0^{r}\kappa (\tau)\mathrm{d}\tau \geq \begin{cases}
       0,  &  2r < \Delta t \\
       r - \Delta t/2, & \mathrm{else}
    \end{cases}
\end{flalign*}
which implies
\begin{flalign}
\label{eq: PiecewiseExp}
    e^{-c\int _0^{r}\kappa (\tau) \mathrm{d}\tau}\leq e^{-c\left(r-\frac{\Delta t}{2}\right)}.
\end{flalign}
In the following we write $x_i$ instead of $x_{i,L}^{u_i}$ and we choose a state feedback control which is given by
\begin{flalign}
    \begin{split}
        u_1(t)&:=(1-\kappa (t))x_L^{0'}(ct)+\kappa (t) x_1(0,t)\\
        u_i(t)&:=-\kappa (t)x_{i-1}'(a_{i-1},t)+\kappa (t) x_i(a_{i-1},t)
    \end{split}
\end{flalign}
where $i\in \{2,\ldots, N_L\}$. Using~\eqref{eq: NeumannSolution} and the variation of constants formula we find, that for $a_{i-1}\leq \omega \leq a_{i-1} + ct$ the solution of the closed-loop system is given by
\begin{flalign}
    \label{eq: SolClosedLoop1}
    \begin{split}
    x_i(\omega,t)&=e^{-c\int _0^{t-\frac{\omega}{c}}\kappa (\tau) \mathrm{d}\tau}x_L^0(a_{i-1})
    +\int _0^{t-\frac{\omega-a_{i-1}}{c}}e^{-c\int _0^{t-\frac{\omega-a_{i-1}}{c}-s}\kappa (\tau) \mathrm{d}\tau}v_i(s)\mathrm{d}s
    \end{split}
\end{flalign}
where
\begin{flalign*}
    v_1(t)\!:= \!(1-\kappa (t))x_L^{0'}(ct),\,\, v_i(t)\!:=\! (1-\kappa (t))x_{i-1}'(a_{i-1},t).
\end{flalign*}
For $a_{i-1} + ct \leq \omega \leq a_i $ the solution is given by
\begin{equation}
    \label{eq: SolClosedLoop2}
    x_i(\omega ,t) := x_L^0(\omega - ct).
\end{equation}
Our aim is, to find an estimate of the form
\begin{equation}
\label{eq: PiecewiseStabEstimate}
    \n{x_i(\cdot,t)}_{H^1(a_{i-1},a_i)}\leq Me^{-kt}\n{x_L^0}_{H^1(a_{i-1},a_i)}
\end{equation}
for some constants $M,k>0$ and for all $i \in\{1,\ldots , N_L\}$. We will proceed by first showing such an estimate for $i=1$.
Using the triangle inequality we find for $a \in [0,a_1]$
\begin{flalign}
\label{eq: Triangle1}
    &\n{x_1(\cdot,t)}_{L^2(0,a)}^2\overset{\eqref{eq: SolClosedLoop1}}{\leq} I_a(t) + J_a(t)
\end{flalign}
where
\begin{flalign}
\label{eq: Triangle2}
\begin{split}
    I_a(t):=\n{e^{-c\int _0^{t-\frac{\cdot}{c}}\kappa (\tau) \mathrm{d}\tau}x_L^0(0)}_{L^2(0,a)}^2, \quad J_a(t):=\n{\int _0^{t-\frac{\cdot}{c}}e^{-c\int _0^{t-\frac{\cdot}{c}-s}\kappa (\tau) \mathrm{d}\tau}v_1(s)\mathrm{d}s}_{L^2(0,a)}^2.
\end{split}
\end{flalign}
These two terms can be estimated via
\begin{flalign}
\label{eq: IEstimate}
\begin{split}
    I_a(t)&=\int _0^{a} e^{-2c\int _0^{t-\frac{\omega}{c}}\kappa (\tau) \mathrm{d}\tau}\n{x_L^0(0)}^2\mathrm{d}\omega 
    \overset{~\eqref{eq: PiecewiseExp}}{\leq} \int _0^{a} e^{-2c\left(t-\frac{\omega}{c}-\frac{\Delta t}{2}\right)}\n{x_L^0(0)}^2\mathrm{d}\omega\\
    &\overset{\omega \leq a}{\leq} \int _0^{a} e^{-2c\left(t-\frac{a}{c}-\frac{\Delta t}{2}\right)}\mathrm{d}\omega \n{x_L^0(0)}^2
    = a e^{2a+c\Delta t} e^{-2ct}\n{x_L^0(0)}^2
\end{split}
\end{flalign}
and
\begin{flalign}
\label{eq: JEstimate}
\begin{split}
J_a(t)
    \!\!=\!\!& \int _0^{a}\n{\int _0^{t-\frac{\omega}{c}}e^{-c\int _0^{t-\frac{\omega}{c}-s}\kappa (\tau) \mathrm{d}\tau}v_1(s)\mathrm{d}s}^2\mathrm{d}\omega
    \!\!\leq \!\! \int _0^{a}\!\!\! e^{-2c\int _0^{t-\frac{\omega}{c}}\kappa (\tau) \mathrm{d}\tau}\int _0^{\frac{\Delta t}{2}}\!\!\!e^{2c\int _{t-\frac{\omega}{c}-s}^{t-\frac{\omega}{c}}\kappa (\tau) \mathrm{d}\tau}v_1(s)^2\mathrm{d}s\mathrm{d}\omega\\
    \leq & \int _0^{a}e^{-2c\left(t-\frac{\omega}{c}-\Delta t\right)}\mathrm{d}\omega \frac{1}{c}\int _0^{\frac{\Delta t}{2}c}x_L^{0'}(s)^2\mathrm{d}s
    \leq  \frac{a}{c}e^{2\left(a+c\Delta t\right)}e^{-2ct}\n{x_L^0}_{H^1(0,a_1)}^2
\end{split}
\end{flalign}
where we used that $v_1(s)=0$ for all $s\geq \frac{\Delta t}{2}$ and $0\leq \kappa (s) \leq 1$ for all $s \geq 0$.
In the following we partition the $H^1$-norm via
\begin{flalign}
\label{eq: NormPartitioning}
\begin{split}
    \n{x_1(\cdot,t)}_{H^1(0,a_1)}^2
    = &\n{x_1(\cdot,t)}_{L^2(0,ct)}^2 + \n{x_1'(\cdot,t)}_{L^2(0,ct)}^2 \quad+\n{x_1(\cdot,t)}_{H^1(ct,a_1)}^2.
\end{split}
\end{flalign}
Also note that due to boundedness of the trace operator on $H^1$ there exists a constant $c_0>0$ such that
\begin{equation}
\label{eq: InitialValueEstimate}
    \forall a\in [\Delta, L_0]\, \forall v\in H^1(0,a): \n{v(0)}^2\leq c_0 \n{v}_{H^1(0,a)}^2.
\end{equation}

We first consider the case $ct \leq a_1$. Applying~\eqref{eq: Triangle1},~\eqref{eq: IEstimate} and~\eqref{eq: JEstimate} to $a=ct$ and using~\eqref{eq: InitialValueEstimate} we find
\begin{flalign*}
    \n{x_1(\cdot,t)}_{L^2(0,ct)}^2 \leq \left( c_0L_0e^{c\Delta t} +\frac{L_0}{c}e^{2c\Delta t}\right)\n{ x_L^{0}}_{H^1(0,ct)}^2.
\end{flalign*}
The derivative of $x_1$ on the interval $[0,ct]$ is given by
\begin{flalign}
\label{eq: DerFormula}
    x_1'(\omega,t)= \kappa\left(t-\frac{\omega}{c}\right)e^{-c\int_0^{t-\frac{\omega}{c}}\kappa (\tau)\mathrm{d}\tau}x_L^0(0)
    \nonumber
    &+\int _0^{t-\frac{\cdot}{c}}\kappa \left(t-\frac
    {\omega}{c}-s\right)e^{-c\int _0^{t-\frac{\cdot}{c}-s}\kappa (\tau) \mathrm{d}\tau}v_1(s)\mathrm{d}s\\
    -&\frac{1}{c}\left(1-\kappa\left(t-\frac{\omega}{c}\right)\right)x_L^{0'}(ct-\omega).
\end{flalign}
The only difference between the first two terms in this formula and $x_1(\omega,t)$ is the scaling term $\kappa$ which is bounded by $0$ and $1$. This observation leads to the estimate
\begin{flalign*}
    \n{x_1'(\cdot,t)}_{L^2(0,ct)}^2 \leq \n{x_1(\cdot,t)}_{L^2(0,ct)}^2+\frac{1}{c}\n{x_L^0(ct-\cdot)}_{L^2(0,ct)}^2
    =\n{x_1(\cdot,t)}_{L^2(0,ct)}^2+\frac{1}{c}\n{x_L^0}_{L^2(0,ct)}^2.
\end{flalign*}
Finally we find
\begin{flalign*}
    &\n{x_1(\cdot,t)}_{H^1(ct,a_1)}^2 = \n{x_L^0(\cdot-ct)}_{H^1(ct,a_1)}^2
    \leq \n{x_L^0}_{H^1(0,a_1)}^2.
\end{flalign*}
Overall we have shown the estimate
\begin{equation*}
    \n{x_1(\cdot,t)}_{H^1(0,a_1)}^2\leq \underset{=:M_1^2}{\underbrace{K_1e^{2a_1}}}e^{-2ct}\n{x_L^0}_{H^1(0,a_1)}^2
\end{equation*}
for $ct \leq a_1$ where
\begin{equation*}
K_1:= 2c_0L_0e^{c\Delta t} +2\frac{L_0}{c}e^{2c\Delta t}+1+\frac{1}{c}.   
\end{equation*}
Now we consider the case $ct > a_1$: The solution $x_1(\cdot,t)$ is given by~\eqref{eq: SolClosedLoop1} on the whole domain $[0,a_1]$.
Applying~\eqref{eq: Triangle1},~\eqref{eq: IEstimate} and~\eqref{eq: JEstimate} to $a=a_1$ and using~\eqref{eq: InitialValueEstimate} we find
\begin{flalign*}
    \n{ x_1(\cdot,t)}_{L^2(0,a_1)}^2 \!\!\leq\!\! \underset{=:K_2}{\underbrace{L_0 (c_0+\frac{1}{c}) e^{2\left(L_0+c\Delta t\right)}}}e^{-2ct}\!\n{ x_L^{0}}_{H^1(0,a_1)}^2.
\end{flalign*}
To estimate the derivative, note that
\begin{flalign*}
    &\n{\left(1-\kappa\left(t-\frac{\cdot}{c}\right)\right)x_L^{0'}(ct-\cdot)}_{L^2(0,a_1)}^2
    =\int_0^{a_1}\left(1-\kappa\left(t-\frac{\omega}{c}\right)\right)^2x_L^{0'}(ct-\omega)^2\mathrm{d}\omega\\
    \leq & \begin{cases}
       \| x_L^{0'} \|_{L^2(0,a_1)}^2, &  t \leq a_1/c + \Delta t / 2 \\
       0, &\textrm{ else}
    \end{cases}
    \leq  e^{2a_1+c\Delta t}e^{-2ct} \| x_L^{0'} \|_{L^2(0,a_1)}^2.
\end{flalign*}
Using~\eqref{eq: DerFormula} we find for the derivative
\begin{equation*}
    \n{x_1'}_{L^2(0,a_1)}^2 \leq \left(K_2+e^{2\left(L_0+c\Delta t\right)}\right)e^{-2ct}\n{x_L^0}_{H^1(0,a_1)}^2.
\end{equation*}
Overall for $ct \geq a_1$ we find the estimate
\begin{flalign*}
    \n{x_1}_{H^1(0,a_1)}^2 \leq \underset{=:M_2^2}{\underbrace{(2K_2+e^{2\left(L_0+c\Delta t\right)})}}e^{-ct}\n{x_L^0}_{H^1(0,a_1)}^2.
\end{flalign*}
For arbitrary $t\geq 0$ we conclude
\begin{equation*}
    \n{x_1}_{H^1(0,a_1)}\leq M e^{-kt}\n{x_L^0}_{H^1(0,a_1)}
\end{equation*}
where $M:= \max \{M_1, M_2\}$ and $k=c$. By considering the solution formulas in~\eqref{eq: SolClosedLoop1} and~\eqref{eq: SolClosedLoop2} we observe, that for $i>1$ the closed-loop solution $x_i$ only depends on 
$x^L_{0_{|{[a_{i-2},a_{i-1}]}}}$. 
Therefore by using analogue estimates as in the case $i=1$ we find that~\eqref{eq: PiecewiseStabEstimate} also holds true for $i>1$. Note that the constants $M$ and $k$ do not depend on the domain length $L$ and the transport equations index $i$. The claim now follows from
\begin{flalign*}
    &\n{x}_{H^1(\Omega ^L)}^2=\sum\nolimits_{i=1}^{N_L}\n{x_i}_{H^1(\Omega _i)}^2
    \leq \sum\nolimits_{i=1}^{N_L}M^2 e^{-2kt}\n{x_0^L}_{H^1(\Omega _i)}^2 = M^2 e^{-2kt} \n{x_0^L}_{H^1(\Omega ^L)}^2.
\end{flalign*}
\hfill
\end{proof}

\section{Numerical example}
\label{Sec: Numerics}
In~\cite{Goettlich2025} it was shown that for a large class of evolution equations domain-uniform stabilizability/detectability implies that exponentially localized perturbations only cause an exponentially localized deviation in the solution of linear quadratic optimal control problems. For the case of Dirichlet boundary control we now present numerical simulations which illustrate that the findings of Section~\ref{Sec: Dirichlet} ensure the same behavior for boundary control systems. To achieve this we solve the optimal control problem
\begin{flalign*}
\begin{split}
    \underset{(x,u)}{\min }\,\,\, \frac{1}{2} \int _0^T \n{x(\cdot,t)}_{L^2(0,L)}^2 \, + \, \alpha &\n{u(\cdot, t)}_{L^2(\Omega _c^L)}^2 \, \mathrm{d}t\\
    \textrm{s.t.}: \, \forall \omega \in \Omega _i \cap \Omega ^L \, \forall t\in [0,T]: \dot{x_i}(\omega,t)&=-cx'(\omega,t)\\
    \forall t \in [0,T]: x_i(a_{i-1},t)-x_{i-1}(a_{i-1},t) &= u_i(t) \\
    \forall \omega \in [0,L]: x_i(\omega,0) &= x_L^0(\omega)|_{\Omega _i \cap \Omega ^L}.
\end{split}
\end{flalign*}
We solve the optimal control problem for two different types of control domains assuming the transport velocity to be constant ($c=2$) in all simulations:
\begin{itemize}
    \item[(1)] For given domain size $L>0$ there are only two control access points at $a_0=0$ and $a_1=\frac{L}{2}$.
    \item[(2)] The control domain consists of a sequence of equidistantly distributed control access points which are given by $a_i=i$ for all $i \in \N_0$.
\end{itemize}

\noindent For the simulation, we use a finite difference method with symmetric difference quotient $(D_hx)(\omega):=\frac{x(\omega + h) - x(\omega -h)}{2h}$ for spatial discretization and an implicit midpoint rule for time discretization of the corresponding optimality system, see, e.g.~\cite[Eq. 5.4]{Goettlich2024}.

\begin{figure}[htb]
        \centering
        \includegraphics[width=\linewidth]{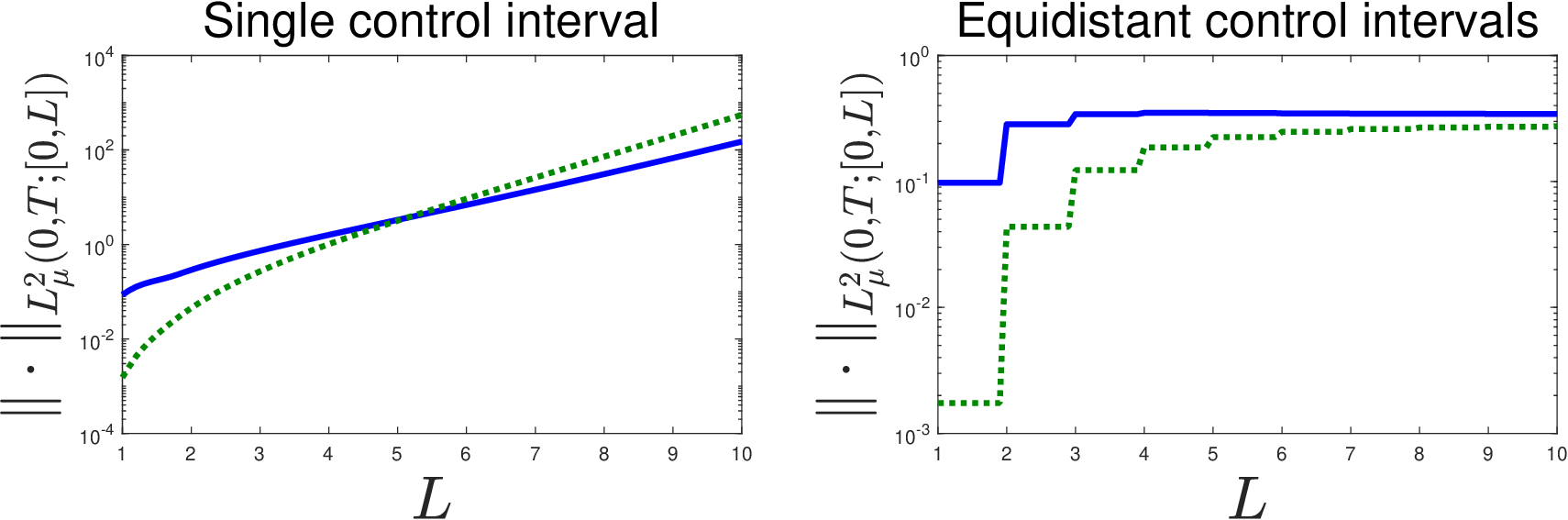}
        \caption{Relation between $L^2(0,T;L_\mu^2(0,L))$-norm of optimal state (solid) and costate (dotted)  and domain size for $T = 5$, $\alpha = 0.156$, $\mu = 0.5$.}
        \label{fig: L2Norm}
\end{figure}

In the following we consider a non-zero initial value as a perturbation of the zero solution of the optimal control problem. Thus, setting parameters $\varepsilon _1 = 0.6$, $\varepsilon _2 = 0.8$ and $\mu(\omega) = 1+\frac{1}{\left(\frac{2}{\varepsilon_2}\left(\omega\!-\!\varepsilon _1\right)\right)^2\!-\!1}$ we choose the initial value
\begin{equation}
\label{Perturbation}
    \varepsilon (\omega):= \begin{cases}
    e^{\mu(\omega)}, & \omega 
    \!\in \!\left (
    \varepsilon_1 - \varepsilon_2/2, \varepsilon_1 + \varepsilon_2/2
    \right )
    \\
    0, & \mathrm{else}
    \end{cases}.
\end{equation}
\begin{figure}[htb]
        \centering
        \includegraphics[trim={0.5cm 0.3cm 1cm 1.3cm},clip,width=\linewidth]{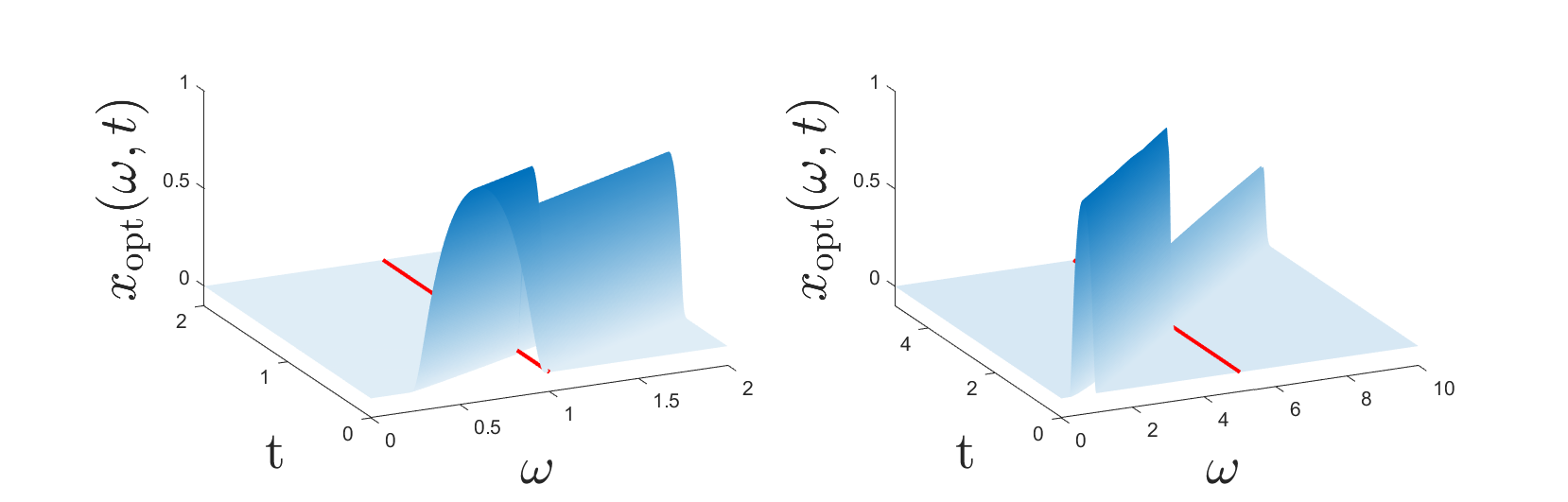}
        \includegraphics[trim={0.5cm 0cm 1cm 1.3cm},clip,width=\linewidth]{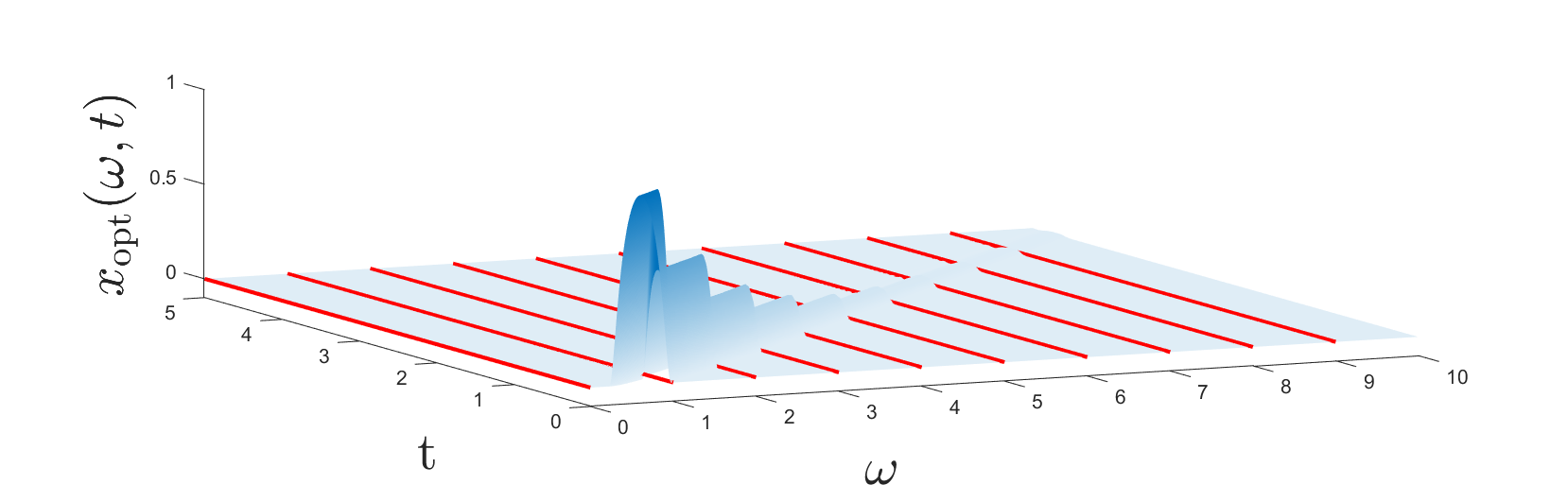}
        \caption{Optimal state trajectory of transport equation with control domain (1) respectively (2) and domain size $L=2$ respectively $L=10$ ($\alpha = 0.156$, $T=5$)}
        \label{fig: 3DPlotState}
\end{figure}
\begin{figure}[htb]
        \centering        \includegraphics[width=\linewidth]{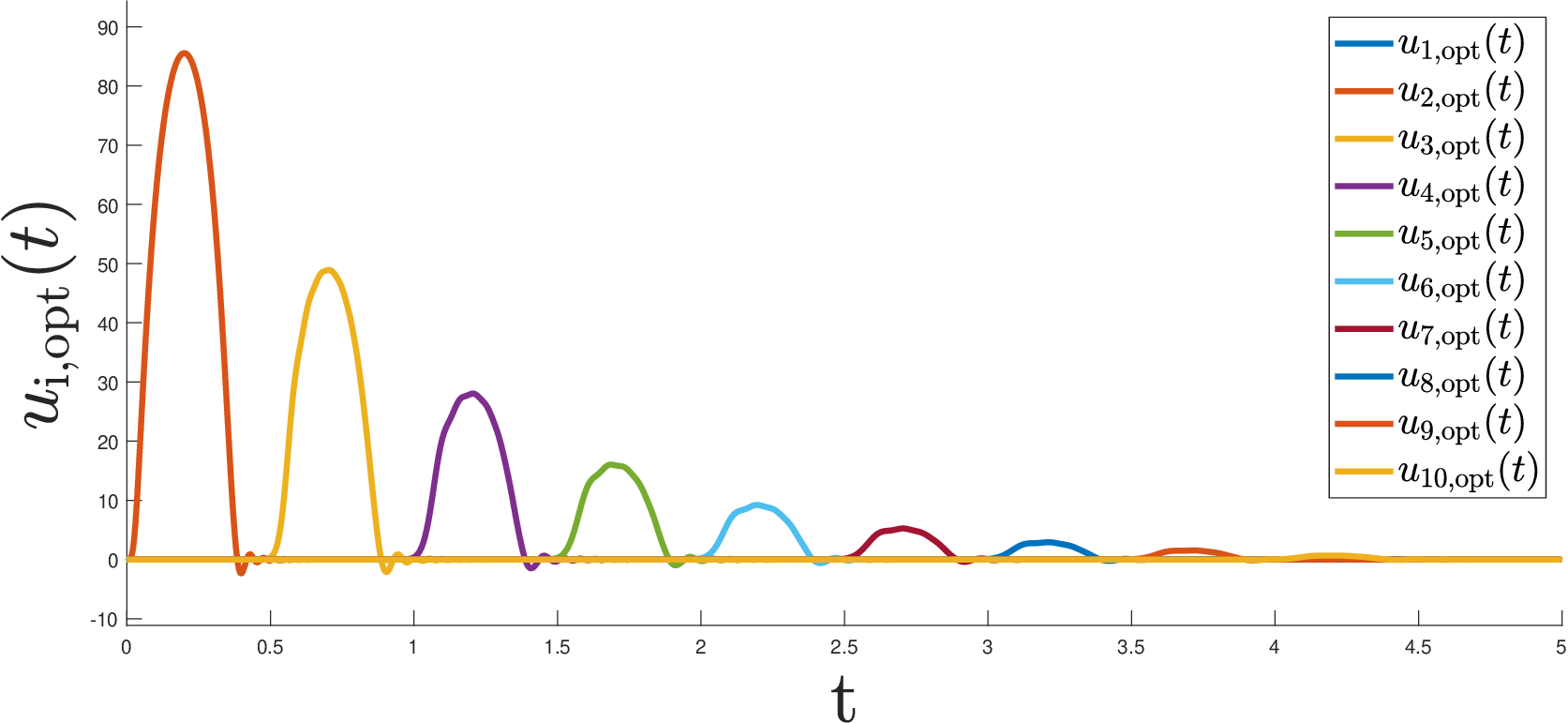}
        \caption{Optimal Dirichlet boundary control input for transport equation with equidistant control access points ($L=10$, $\alpha = 0.156$, $T=5$)}
        \label{fig: OptInput}
\end{figure}
Figure~\ref{fig: 3DPlotState} shows, how the dynamics of the system transport the initial value perturbation along the whole length of the spatial domain. At each control access point part of the perturbation is compensated which leads to an exponential decay in space. However the optimal control does not directly eliminate the perturbation once it reaches the control domain (see also Figure~\ref{fig: OptInput}). This is due to the strictly positive weight $\alpha$ which introduces an input cost. 
The upper left plot in Figure~\ref{fig: 3DPlotState} shows the optimal state trajectory on the domain $L=2$, where the control domain choices of (1) and (2) coincide. For the upper right and bottom plot the domain size was increased to $L=10$. In the case of a single control access point (upper right) the exponential decay becomes a lot slower while in the case of equidistant control access points it the decay rate remains unchanged (bottom).
In Figure~\ref{fig: L2Norm} we illustrated this behaviour more explicitly by plotting the (spatially) exponentially weighted $L^2$-norm
\begin{flalign*}
    \| x \|_{L^2(0,T;L_\mu^2(0,L))} & := \| e^{\mu\n{\omega - \varepsilon _1}_1} x(t,\omega) \|_{L^2(0,T;L_\mu^2(0,L))}
\end{flalign*}
for increasing domain sizes $L$. In the case of a single control access point this quantity increases exponentially while in the case of equidistant control access points it reaches an upper bound. This illustrates the major impact of domain-uniform stabilizability on the exponential decay of the influence of exponentially localized perturbations on the optimal solutions behaviour.

\section{Conclusions}
In this paper we derived a simple characterization of domain-uniform stabilizability for a chain of transport equations with a Dirichlet/Neumann type control. 
In future work more complex networks of transport equations could be considered.

\section*{Acknowledgements}
We thank Timo Reis and Thavamani Govindaraj (TU Ilmenau) for fruitful discussions on a suitable formulation of the transport equation with Neumann boundary control.

\bibliographystyle{abbrv}
\bibliography{References}
\end{document}